\documentclass[default]{sn-jnl}

\usepackage{amsmath}
\usepackage{graphicx}%
\usepackage{multirow}%
\usepackage{amsmath,amssymb,amsfonts}%
\usepackage{amsthm}%
\usepackage{mathrsfs}%
\usepackage[title]{appendix}%
\usepackage{xcolor}%
\usepackage{textcomp}%
\usepackage{manyfoot}%
\usepackage{booktabs}%
\usepackage{algorithm}%
\usepackage{algorithmicx}%
\usepackage{algpseudocode}%
\usepackage{listings}%

\theoremstyle{thmstyleone}%
\newtheorem{theorem}{Theorem}
\newtheorem{proposition}[theorem]{Proposition}%

\theoremstyle{thmstyletwo}%

\theoremstyle{thmstylethree}%
\newtheorem{definition}{Definition}%

\newtheorem{lemma}[theorem]{Lemma}

\begin{document}

	\title[Existence of weak solutions for fractional $p(x,.)$ -Laplacian Dirichlet problems with nonhomogeneous boundary conditions.]{Existence of weak solutions for fractional $p(x,.)$-Laplacian Dirichlet problems with nonhomogeneous boundary conditions.}


	\author*[1]{\fnm{Achraf} \sur{El wazna}}\email{elwaznaachraf@gmail.com}

\author[1]{\fnm{Azeddine} \sur{Baalal}}\email{abaalal@gmail.com}

\affil*[1]{\orgdiv{Department of Mathematics and Computer Science}, \orgname{University Hassan II, Faculty of Sciences A\"{i}n Chock}, \orgaddress{\street{Road EL Jadida Km 8}, \city{Casablanca}, \postcode{20100}, \state{B.P. 5366 M\^{a}arif}, \country{Morocco}}}


\abstract{In this paper, we consider the existence  of solutions of the following nonhomogeneous fractional  $p(x,.)$-Laplacian Dirichlet problem:
	\begin{equation*}
	\left\{\begin{aligned}
	\Big(-\Delta_{p(x,.)}\Big)^s u (x)&=f(x, u)  &\text { in }& \Omega, \\
	u &=g &\text { in }& \mathbb{R}^N \setminus\Omega,
	\end{aligned}\right.
	\end{equation*}
	where $\Omega\subset\mathbb{R}^N$ is a smooth bounded domain, $\Big(-\Delta_{p(x,.)}\Big)^s$ is the fractional $p(x,.)$-Laplacian, $f$ is a Carathéodory function with suitable growth condition and $g$ is a given boundary data. The proof of our main existence results relies on the study of the fractional $p(x, \cdot)$-Poisson equation with a nonhomogeneous Dirichlet boundary condition and the theory of fractional Sobolev spaces with variable exponents, together with Schauder's fixed point theorem.}




\maketitle
\section{Introduction}

In this paper, we are interested in the existence of weak solutions for the following problem:
\begin{equation*}\label{P}\tag{P}
\left\{\begin{aligned}
\Big(-\Delta_{p(x,.)}\Big)^s u (x)&=f(x, u)  &\text { in }& \Omega, \\
u &=g &\text { in }& \mathbb{R}^N \setminus\Omega,
\end{aligned}\right.
\end{equation*}
where $\Omega\subset \mathbb{R}^N$ is a smooth bounded domain, $s\in (0,1)$,   $f:\Omega\times\mathbb{R}\mapsto \mathbb{R}$ is a  Carathéodory function, $g$ is a given boundary data, and $\Big(-\Delta_{p(x,.)}\Big)^s$ is the fractional $p(x,.)$-Laplacian operator defined as 
\begin{equation}\label{fractionallaplacian}
\Big(-\Delta_{p(x,.)}\Big)^s u(x)=p.v.\int_{\mathbb{R}^N} \dfrac{\lvert u(x)-u(y)\lvert ^{p(x,y)-2}(u(x)-u(y))}{\lvert x-y\lvert ^{N+sp(x,y)}} dy,
\end{equation}
for all $x\in\mathbb{R}^N$, where $p:\mathbb{R}^N\times \mathbb{R}^N\mapsto (1,\infty)$ is bounded, continuous and symmetric function.

To the best of our knowledge, this operator was first introduced by Kaufmann, Rossi and Vidal \cite{zbMATH06931307}, in which the authors extended the Sobolev spaces with variable exponents to include the fractional case together with a compact embedding theorem. As an application, they established the existence and uniqueness of a solution for the following fractional $p(x,.)$-Laplacian problem using a direct method of calculus of variations.:
\begin{equation*}
\left\{\begin{aligned}
\Big(-\Delta_{p(x,.)}\Big)^s u (x)+ |u(x)|^{q(x)-2}u(x)&=f(x)  &\text { in }& \Omega, \\
u &=0 &\text { on }& \partial\Omega,
\end{aligned}\right.
\end{equation*} 
where $q(x)>p(x,x)$ and $f\in L^{a(x)}(\Omega)$ for some $a(x)>1$ .\\

In \cite{zbMATH06810426}, Bahrouni and R{\u{a}}dulescu obtained some qualitative properties both on the fractional Sobolev space with variable exponent $W^{s,q(.),p(.,.)}(\Omega)$ and the related nonlocal operator $\Big(-\Delta_{p(x,.)}\Big)^s$ . Moreover, by using Ekeland's variational method, they studied the existence of solutions for the following problem:
\begin{equation*}
\left\{\begin{aligned}
\Big(-\Delta_{p(x,.)}\Big)^s u (x)+ |u(x)|^{q(x)-1}u(x)&=\lambda |u(x)|^{r(x)-1}u(x)  &\text { in }& \Omega, \\
u &=0 &\text { on }& \partial\Omega,
\end{aligned}\right.
\end{equation*} 
where $\lambda>0$, and $\displaystyle 1<r(x)<p^{-}:=\min_{(x,y)\in\Omega\times\Omega}p(x,y)$.\\

Bahrouni \cite{zbMATH06813457} continued the study of this class of fractional Sobolev spaces with variable exponent and the fractional $p(x,.)-$Laplacian. More precisely, he showed a variant of comparison principle for the operator $\Big(-\Delta_{p(x,.)}\Big)^s$, and proved a general principle of sub-supersolution method for problem \eqref{P} with homogeneous boundary condition, i.e., $g=0$, where $p,f$ are continuous functions and $f$ satisfies the following assumption:
\begin{equation*}
|f(x,t)|\leqslant c_1+c_2|t|^{r(x)-1},\hspace{1cm} \forall x\in \mathbb{R}^N,\ \forall t\in\mathbb{R},
\end{equation*}
where $r\in C(\mathbb{R}^N,\mathbb{R})$ and $1<r(x)<\bar{p}^*_s:=\frac{Np(x,x)}{N-sp(x,x)}\hspace{0.5cm}\forall x\in \mathbb{R}^N.$\\

Subsequently, Ho and Kim \cite{zbMATH07134940} refined the fractional Sobolev space with variable exponent given in \cite{zbMATH06931307} and established fundamental embeddings for these new spaces, which generalize the well-known fractional Sobolev spaces $W^{s,p}(\Omega)$. With these embeddings, and by providing some conditions for the exponent  $p(.,.)$ in $\mathbb{R}^N\times\mathbb{R}^N$, they investigate the boundedness and the multiplicity of the solutions to the problem \eqref{P}, in the case of  boundary data $"g=0"$, using De Giorgi's iteration and the localization method. Afterwards, the  fractional Sobolev spaces with variable exponents and problems involving the operator $(-\Delta_{p(x,.)})^s$ has been extensively studied in various contexts. Readers may refer to \cite{vplapp1,vplapp2,vplapp3,vplapp4,vplapp5,vplapp6} and the references therein.\\

When $p(x,y)\equiv p$ is constant, the operator \eqref{fractionallaplacian} reduces to the fractional $p$-Laplacian $\big(-\Delta\big)^s_p$, defined by 
\begin{equation*}
	(-\Delta)^s_p u(x)=p.v.\int_{\mathbb{R}^N} \dfrac{\lvert u(x)-u(y)\lvert ^{p-2}(u(x)-u(y))}{\lvert x-y\lvert ^{N+sp}} dy.
\end{equation*}
This operator is consistent, up to a normalization constant depending on $N$ and $s$, with the linear fractional Laplacian $\big(-\Delta\big)^s$ in the case $p=2$.
This kind of nonlocal operator has a wide range of real-world applications. It naturally arises in several physical phenomena, such as flames propagation and chemical reactions of liquids, in population dynamics and geophysical fluid dynamics, or in mathematical finance. For an introduction to this topic and further references on nonlocal operators and their applications, see \cite{hitchh} and the references therein.\\

Note that in the constant case the problem \eqref{P} becomes 
\begin{equation*}\label{Pconstant}\tag{$P_1$}
\left\{\begin{aligned}
(-\Delta)^s_p u (x)&=f(x, u)  &\text { in }& \Omega, \\
u &=g &\text { in }& \mathbb{R}^N \setminus\Omega.
\end{aligned}\right.
\end{equation*}
A broad range of existence results for the problem \eqref{Pconstant}, with homogeneous Dirichlet boundary conditions (i.e. $g=0$), has been obtained by Iannizzotto et al. in \cite{zbMATH06567151} via Morse theory under different growth assumptions for $f$. We can refer to \cite{plapp1,plapp2,plapp3,plapp4,plapp5,plapp6,plapp7} for existence, multiplicity and regularity results for fractional $p$-Laplacian problems. For the special case $p=2$, several existence results via variational methods are proved in a series of papers by Servadei and Valdinocci \cite{serval3,serval2,serval1,serval4}. See also \cite{lapp1,lapp2,lapp3,lapp4} for the related discussions.\\

The operator $(-\Delta_{p(x,.)})^s$ can be seen as the fractional version of the well-known $p(x)-$Laplacian, given by $\Delta_{p(.)}u=div(|\nabla u|^{p(.)-2}\nabla u)$ associated with the variable exponent Sobolev space $W^{1,p(.)}(\Omega)$. These spaces have attracted increasing interest in the last two decades as they provide a natural setting for studying differential equations and variational problems with variable exponent growth conditions. 
The study of these problems is important not only for mathematical purposes but also for their significance in real-world applications, such as image processing, nonlinear electrorheological fluids, and elastic mechanics (see the monograph \cite{zbMATH06435536} and the references therein).\\
Note that the problem \eqref{P} can be seen as the fractional form of the following problem: \begin{equation*}\label{Pvar}\tag{$P_2$}
\left\{\begin{aligned}
\Big(-\Delta _{p(x)}\Big) u (x)&=f(x, u)  &\text { in }& \Omega, \\
u &=g &\text { on }& \partial\Omega.
\end{aligned}\right.
\end{equation*}
There is an extensive literature dealing with the problem \eqref{Pvar}, mostly with homogeneous boundary data. For results concerning the existence of solutions to these types of problems, we refer the readers to  \cite{pxlapp1,pxlapp2,pxlapp3,pxlapp4,pxlapp5,pxlapp6,zbMATH01891301} (see \cite{defigue} for the constant exponent case $p(.)\equiv p$).\\
Let us also point out that in \cite{MR3911965}, using an approximation argument, Baalal and Berghout prove that for every continuous function $g$ on the boundary $\partial\Omega$ there exists a unique continuous solution of the problem \eqref{Pvar}. Notice that our  situation is different from \cite{MR3911965}, because of the presence of a nonlocal  operator, in the sens that the value of $(-\Delta_{p(x,.)})^s u(x)$ at any point $x\in\Omega$ depends not only on the value of $u$ on $\Omega$, but actually on its value on all of $\mathbb{R}^N$.\\

Our interest in Problem \eqref{P} has originated from the paper \cite{MR1837267}, in which Baalal proves the existence and uniqueness of a weak solution to the Dirichlet problem  
\begin{equation}\label{baalal}\tag{$P_3$}
\left\{\begin{aligned}
\mathcal{L}u:= \displaystyle\sum_j \frac{\partial}{\partial x_j}\Big(\sum_i\alpha_{i,j}\frac{\partial u}{\partial x_i}+\delta_i u\Big)&=f(.,u,\nabla u) \quad \text { in } \Omega, \\
u &=g \text { on } \partial \Omega  .
\end{aligned}\right.
\end{equation}
Through the change of unknown $v=u-\tilde{g}$, where $\tilde{g}$ is an extension of $g$ to $\Omega$, the nonhomogeneous Dirichlet problem given in \eqref{baalal} reduced into a homogeneous Dirichlet problem. However the nonlinear case is quite different,
in our case, an analogous change of unknown function leads to consider the following operator $v\longrightarrow \Big(-\Delta_{p(x,.)}\Big)^s(v+\tilde{g})$. This operator turns out to be harder to handle than the fractional $p(x,.)$-Laplacian, and for that matter, we will deal directly with the nonhomogeneous Dirichlet Problem \eqref{P}.\\

The proof of our existence result relies on a combination of the Schauder's fixed point theorem, in the spirit of \cite{MR1837267}, and an $L^{r(\cdot)}$-estimate of the solution to the following fractional nonhomogeneous Poisson equation associated with Problem \eqref{P}:
\begin{equation}\label{poi}\tag{$P_{h,g}$}
\left\{\begin{aligned}
\Big(-\Delta_{p(x,.)}\Big)^s u (x)&=h(x)  &\text { in }& \Omega, \\
u &=g &\text { in }& \mathbb{R}^N \setminus\Omega,
\end{aligned}\right.
\end{equation}
We believe that this result has independent value and can be applied in various other settings. Specifically, Section \ref{sec3} can be considered as a self-contained and independent part of our work.\\

Having at hand a regularity result for the problem \eqref{poi} and inspired by \cite{MR1837267}, we develop a fixed point argument to obtain a solution of \eqref{P}. Note that, due to the nonlocality and nonlinearity of the operator and the nonhomogeneous Dirichlet condition, the approach of \cite{MR1837267} has to be adapted significantly. The main feature is the consideration of the fractional $p(x,.)$-Laplacian together with a nonzero boundary condition.\\ 

This paper is organized as follow. In section \ref{sec2}, we briefly describe the natural functional framework for our problem including the fractional Sobolev spaces with variable order. In Section \ref{sec3} we study the fractional $p(x,.)$-Poisson's Problem \eqref{poi}, with nonhomogeneous Dirichlet condition. Finally in Section \ref{sec4}, we establish the main existence result for Problem  \eqref{P}.
\section{Preliminary Definitions and Results}\label{sec2}
In this section, we recall some definitions and basic properties of Lebesgue and fractional Sobolev spaces with variable exponents. We also provide some necessary lemmas that will be used in this paper.\\
\subsection{Variable exponent Lebesgue spaces}
Let $\Omega$ be a smooth bounded domain in $\mathbb{R}^{N}$, we consider the set

$$C_{+}(\overline{ \Omega})=\left\{h\in C(\overline{ \Omega}) : 1<\underset{x\in\overline{\Omega}}{\inf}~~ h(x)\leqslant \underset{x\in\overline{\Omega}}{\sup}~~ h(x)<\infty  \right\}.$$ 
For any $ h\in C_{+}(\overline{ \Omega})$, we denote 
$$h^{+}= \underset{x\in\overline{\Omega}}{\sup} ~~h(x)~~ \text{ and}~~ h^{-}=\underset{x\in\overline{\Omega}}{\inf}~~ h(x).$$
For $q\in C_{+}(\overline{ \Omega})$ , we define the variable exponent Lebesgue space $L^{q(.)}(\Omega) $ as 

$$ L^{q(.)}(\Omega)= \left\{u: \Omega \to \mathbb{R}~\text{is measurable with } \int_{\Omega}\lvert u(x)\lvert^{q(x)}dx < \infty  \right\},$$
which is endowed with the so-called Luxembourg norm
 \begin{equation*} 
\| u\|_{L^{q(.)}(\Omega)}:=\| u\|_{q(.),\Omega}= \inf \left\{\lambda > 0 : \int_{\Omega}\left\lvert\dfrac{u(x)}{\lambda}\right\lvert^{q(x)}dx \leqslant  1  \right\}.
\end{equation*}
It is well known that $\big(L^{q(.)}(\Omega), \|\ .\ \|_{q(.),\Omega}\big)$ is a separable, uniformly convex Banach space (see \cite{zbMATH00095260, zbMATH01703050}).\\
Let $q^{\prime}\in C_{+}(\overline{ \Omega}) $ be the conjugate exponent of $q$, that is, $\frac{1}{q(x)}+\frac{1}{q^{\prime}(x)}=1$. Then we have the following Hölder type inequality.
\begin{lemma}[\cite{zbMATH00095260}]\label{Holder} If $u\in L^{q(.)}(\Omega)$ and $v\in L^{q^{\prime}(.)}(\Omega)$, then
	\begin{equation*}
	\left\lvert\int_{\Omega}u.v dx\right\rvert \leqslant  \left( \dfrac{1}{q^-}+\dfrac{1}{q^{\prime -}}\right)\| u\|_{q(.),\Omega} \| v\|_{q^{\prime}(.),\Omega} \leqslant 2 \| u\|_{q(.),\Omega} \| v\|_{q^{\prime}(.),\Omega}.
	\end{equation*}    
\end{lemma}
\begin{lemma}[\cite{variablemono1} Corollary 2.28]\label{holder2}
 Let $p, q,\beta\in C_{+}(\overline{ \Omega}) $ be such that $$\frac{1}{\beta(x)}=\frac{1}{p(x)}+\frac{1}{q(x)}$$
 $\text{for any}\ x\in\overline{ \Omega}$. Then 
 $$\|uv\|_{\beta(.),\Omega}\leqslant \bigg(\Big(\frac{\beta}{p}\Big)^++\Big(\frac{\beta}{q}\Big)^+\bigg)\|u\|_{p(.),\Omega}\|v\|_{q(.),\Omega},$$
 for all $u\in L^{p(.)}(\Omega)$ and $v\in L^{q(.)}(\Omega).$
\end{lemma}
\vspace{0.3cm}
\begin{lemma}\label{inclusion}
	Let $p_1,p_2\in C_+(\overline{ \Omega})$, such that $p_1(x)\leqslant p_2(x)\ \text{ for all  } x\in \overline{ \Omega},$ then $L^{p_2(.)}(\Omega)$ is continuously embedded in $L^{p_1(.)}(\Omega)$. 
\end{lemma}
\vspace{0.3cm}

\begin{proposition}[\cite{zbMATH01703050}, Theorem 1.16]\label{nemytsky}
	If $f: \Omega\times\mathbb{R}\rightarrow\mathbb{R}$ is a Carathéodory function and satisfies \begin{equation*}
		|f(x,t)|\leqslant a(x)+ b|t|^{p_1(x)/p_2(x)}\ \text{ for all  } x\in\Omega,\ t\in\mathbb{R},
	\end{equation*}
	where $p_1,p_2\in C_+(\Omega),$ $a\in L^{p_2(.)}(\Omega)$, $a(x)\geqslant 0$ and $b\geqslant 0$ is constant, then the Nemytsky operator $\mathcal{N}_f:\ L^{p_1(.)}(\Omega) \rightarrow L^{p_2(.)}(\Omega)$ defined by $\big(\mathcal{N}_fu\big)(x)=f(x,u(x))$, is a continuous and bounded operator.
\end{proposition}

An important role in manipulating the variable Lebesgue spaces is played by the modular of the  $L^{q(.)}(\Omega)$ space, which is the mapping  $\rho_{q(.)}: L^{q(.)}(\Omega) \longrightarrow \mathbb{R}$ defined by 
\begin{equation}
\rho_{q(\cdot)}(u):=\int_{\Omega}\lvert u(x)\lvert ^{q(x)} d x .
\end{equation}
From this, we have the following relations between the norm and modular (see \cite{zbMATH01703050,zbMATH01891301}).
\begin{proposition}  Let $u \in L^{q(.)}(\Omega)$ and $\left\{u_{j}\right\}_{j} \subset L^{q(.)}(\Omega)$, then
	\begin{enumerate}
		\item $\|u\|_{q(.),\Omega}<1$ (resp. $=1,>1$ ) $\Leftrightarrow \rho_{q(.)}(u) < 1$ ( resp. $=1,>1$),
		
		\item $\|u\|_{q(.),\Omega}  > 1 \Rightarrow\|u\|_{q(.),\Omega}^{q^{-}} \leqslant  \rho_{q(.)}(u) \leqslant \|u\|_{q(.),\Omega}^{q+}$,
		
		\item $\|u\|_{q(.),\Omega} <1 \Rightarrow \|u\|^{q+}_{q(.),\Omega} \leqslant  \rho_{q(.)}(u) \leqslant \|u\|_{q(.),\Omega}^{q^{-}}$,
		
		\item $\displaystyle\lim _{j \rightarrow \infty}\|u_j\|_{q(.),\Omega}=0(\infty) \Leftrightarrow \lim _{j \rightarrow \infty} \rho_{q(.)}\left(u_{j}\right)=0(\infty)$,

		\item $\displaystyle\lim _{j \rightarrow \infty} \| u_{j}-u\|_{q(.),\Omega}=0 \Leftrightarrow \lim _{j \rightarrow \infty} \rho_{q(.)}\left(u_{j}-u\right)=0$.
	\end{enumerate}
\end{proposition}
\begin{proposition}
 Let $u \in L^{q(.)}(\Omega)$ and $\left\{u_{j}\right\}_{j} \subset L^{q(.)}(\Omega)$, then the following statements are equivalent each other:\begin{enumerate}
 	\item $\displaystyle\lim _{j \rightarrow \infty} \| u_{j}-u\|_{q(.),\Omega}=0 $.
 	\item $\displaystyle\lim _{j \rightarrow \infty} \rho_{q(.)}\left(u_{j}-u\right)=0$.
 	\item $u_j \to u$ in measure in $\Omega$ and $\displaystyle\lim_{j\to \infty}\rho_{q(.)}\left(u_{j}\right)=\rho_{q(.)}\left(u\right)$
 \end{enumerate}
\end{proposition}
The following lemma will be used in the sequel.
\begin{lemma}[\cite{zbMATH05854634}, Lemma 3.2.12 ]\label{Cara}
	Let $\gamma\in  C_{+}(\overline{ \Omega})$. Then 
	$$\min\big\{\lvert \Omega\lvert ^{\frac{1}{\gamma^+}}; \lvert \Omega\lvert ^{\frac{1}{\gamma^-}}\big\}\leqslant \|1\|_{\gamma(.),\Omega}\leqslant \max\big\{\lvert \Omega\lvert ^{\frac{1}{\gamma^+}}; \lvert \Omega\lvert ^{\frac{1}{\gamma^-}}\big\},$$
	where $|\Omega|$ is the Lebesgue measure of $\Omega$.
\end{lemma}
For the development and properties concerning the variable exponent Lebesgue space $L^{p(.)}(\Omega)$ one may
refer to \cite{variablemono1,zbMATH05854634} and the references therein.
\subsection{The fractional Sobolev spaces with variable exponents}
In this subsection, we recall some useful properties of the fractional Sobolev spaces with variable exponents. We refer the reader to \cite{zbMATH06931307,zbMATH06810426,zbMATH07134940,generalazroul,vplapp3} for more details.\\
Let  $0<s<1$ and let $\mathfrak{C}(\overline{ \Omega})$ the set of functions  $p\in C\big(\overline{ \Omega}\times\overline{ \Omega}\big)$ satisfying 
\begin{equation}\label{expcond1}
1<p^{-}:=\underset{\overline{\Omega} \times \overline{\Omega} }{\inf}~~ p(x,y)\leqslant p^{+}:= \underset{\overline{\Omega} \times \overline{\Omega} }{\sup}~~ p(x,y)<\infty
\end{equation}
and 
\begin{equation}\label{expcond2}
	p\  \text{is symmetric, that is, }\ p(x,y)=p(y,x) \text{for all}\ (x,y)\  \in \overline{\Omega}\times\overline{\Omega}.
\end{equation}\\
For $q\in C_+(\overline{\Omega})$ and $p\in \mathfrak{C}(\overline{ \Omega})$, we define the fractional Sobolev space with variable exponent $W^{s,q(.),p(.,.)}(\Omega)$ as follows: 
\begin{equation*}
W^{s,q(.),p(.,.)}(\Omega):= \left\{u\in L^{q(.)}(\Omega): \rho_{s,p(.,.),\Omega}(u) < \infty \right\},
\end{equation*}
where the modular $\rho_{s,p(.,.),\Omega}$ is defined as 
$$\displaystyle\rho_{s,p(.,.),\Omega}(u)=\int _{\Omega}\int_{\Omega} \dfrac{\lvert u(x)-u(y)\lvert ^{p(x,y)}}{\lvert x-y\lvert ^{N+sp(x,y)}} dxdy.$$
The space $W^{s,q(.),p(.,.)}(\Omega)$ is a separable reflexive Banach space if it is endowed with the norm
\begin{equation}
\|u\|_{s,q,p,\Omega}=[u]_{s,p(.,.),\Omega}+\| u\|_{q(.),\Omega}\  ,
\end{equation}
where $[u]_{s,p(.,.),\Omega}$ is the corresponding variable exponent Gagliardo seminorm which is defined by $$[u]_{s,p(.,.),\Omega}= \inf \left\{\lambda > 0 :  \rho_{s,p(.,.),\Omega}\bigg(\frac{u}{\lambda}\bigg) \leqslant 1 \right\}.$$
The following proposition states the relation between the seminorm $[\ .\ ]_{s,p,\Omega}$ and the modular $\rho_{s,p(.,.)}$.
\begin{proposition} \label{Pro} Let $u\in W^{s,q(.),p(.,.)}(\Omega)$, then 
	\begin{enumerate}
		\item $[u]_{s,p(.,.),\Omega}\geqslant 1 \Rightarrow[u]_{s,p(.,.),\Omega}^{p^{-}} \leqslant \rho_{s,p(.,.),\Omega}(u) \leqslant [u]_{s,p(.,.),\Omega}^{p^{+}}$,
		\item $[u]_{s,p(.,.),\Omega}<1 \Rightarrow [u]_{s,p(.,.),\Omega}^{p^{+}} \leqslant \rho_{s,p(.,.),\Omega} (u)\leqslant [u]_{s,p(.,.),\Omega}^{p^{-}}$.
	\end{enumerate}
\end{proposition}
In the following, to shorten notation, we write $\overline{p}(x)$ instead of $p(x,x)$ and hence $\overline{p}\in C_+(\overline{ \Omega})$. Further if $q(x)=\overline{p}(x) $, we simply write $W^{s,p(.,.)}(\Omega)$ and  $ \|u\|_{s,p(.,.),\Omega}$ instead of $W^{s,q(.),p(.,.)}(\Omega)$  and $  \|u\|_{s,q(.),p(.,.),\Omega}$, respectively.\\

In \cite{zbMATH06931307} the authors proved  a continuous and compact embedding theorem for the spaces $W^{s,q(.),p(.,.)}(\Omega)$ under the assumption $q(x)>\overline{p}(x)$ on $\overline{ \Omega}$. The authors in \cite{zbMATH07134940} give an improvement of \cite[Theorem 1.1]{zbMATH06931307} assuming that $q(x)\geqslant\overline{p}(x),$ for all $x\in\overline{ \Omega}.$\\

\begin{theorem}[\cite{zbMATH07134940}]\label{emb}
	Let $\Omega\subset \mathbb{R}^N$ be a smooth bounded domain, $s\in(0,1)$ , $q\in C_+(\overline{ \Omega})$ and $p\in \mathfrak{C}(\overline{\Omega})$, such that $sp^+<N$ and $q(x)\geqslant p(x,x)$ for all $x\in\overline{ \Omega}$.\\
	Assume that $r \in C_{+}(\overline{\Omega})$,  such that 
	$$ r(x) <p^{*}_{s}(x):= \dfrac{N\bar{p}(x)}{N-s\bar{p}(x)},\ \text{for all}\  x\in \overline{\Omega}.$$
	Then, there exists a constant $C=C(N,s,p,q,r,\Omega)>0$ such that for any $u\in W^{s,q(.),p(.,.)}(\Omega)$ 
	$$ \|u\|_{r(.),\Omega}\leqslant C\| u\| _{s,q,p,\Omega}.$$
	That is, the space $W^{s,q(.),p(.,.)}(\Omega)$ is continuously embedded in $L^{r(x)}(\Omega)$. Moreover, this embedding is compact.
\end{theorem}
 Note that if $q(x)\geqslant \overline{p}(x)$ on $\overline{ \Omega}$, the space $W^{s,q(.),p(.,.)}(\Omega)$ is continuously embedded in $W^{s,p(.,.)}(\Omega)$. \\
 
We denote by $$W^{s,p(.,.)}_{0}(\Omega)=\left\{u\in W^{s,p(.,.)}(\mathbb{R}^N):\ u=0\ on\ \mathbb{R}^N\setminus\Omega \right\}.$$
By the same argument as in Lemma 3.1 in \cite{zbMATH06813457} and using Theorem \ref{emb} we can obtain the following Poincaré type inequality.
\begin{lemma}\label{poincaré}
	Let $\Omega\subset \mathbb{R}^N$ be a smooth bounded domain, $s\in(0,1)$ and $p\in \mathfrak{C}(\mathbb{R}^N)$, with $sp^+<N$.
Then, there exists a constant $C>0$ such that 
$$ \|u\|_{\overline{p}(.),\Omega}\leqslant C [u]_{s,p(.,.),\mathbb{R}^N}.$$
\end{lemma}

The following technical lemmas will be useful in the proof of our main result.
\begin{lemma}[\cite{zbMATH01548957}, Lemma 2.1]\label{edm}
	Let $\alpha,\beta\in C_+(\overline{ \Omega})$ and $u\in L^{\alpha (.)}(\Omega)$. Then 
	\begin{eqnarray*}
	\text{If}&  \|u\|_{\alpha\beta(.),\Omega}\leqslant 1,\ &\text{ then}\hspace{0.5cm}	\|u\|^{\beta^+}_{\alpha\beta(.),\Omega}\leqslant\||u|^{\beta(.)}\|_{\alpha(.),\Omega}\leqslant\|u\|^{\beta^-}_{\alpha\beta(.),\Omega}. \\
			\text{If}&  \|u\|_{\alpha\beta(.),\Omega}\geqslant 1,\ &\text{ then}\hspace{0.5cm}	\hspace{0.5cm}	\|u\|^{\beta^-}_{\alpha\beta(.),\Omega}\leqslant\||u|^{\beta(.)}\|_{\alpha(.),\Omega}\leqslant\|u\|^{\beta^+}_{\alpha\beta(.),\Omega}.
	\end{eqnarray*}
\end{lemma}
\begin{lemma}[\cite{vplapp3}]\label{tracebergot}
	Let $\Omega\subset \mathbb{R}^N$ be an open subset, and $p\in \mathfrak{C}(\mathbb{R}^N)$, such that
	$sp^+<N$.
	Assume that $u\in W^{s,p(.,.)}(\Omega)$. If there exists a compact subset $K\subset \Omega$ such that $u= 0 \text{ in } \Omega \setminus K $, then the extension function $\widetilde{u_\Omega}$ defined as  
	
	\begin{equation*}
	\widetilde{u_\Omega}(x)=\left\{\begin{aligned}
	& u(x)\ &\text{if} ~~& x\in \Omega,\\
	&~0\  &\text{if} ~~& x\in \mathbb{R}^N\setminus\Omega\\
	\end{aligned}\right.
	\end{equation*}
	belongs to $W^{s,p(.,.)}_0(\Omega)$.
\end{lemma}
\begin{lemma}\label{dx}
	Let $\Omega$ a bounded domain,  and 
$p\in \mathfrak{C}(\mathbb{R}^N)$, with $sp^+<N$. Assume that  $\varphi\in W^{s,p(.,.)}(\Omega)$. Then there exists a sequence $\{\varphi_n\}_{n\in\mathbb{N}}\subset W^{s,p(.,.)}_0(\Omega)$ such that 
	$$\lvert \varphi_n(x)\lvert \nearrow \lvert \widetilde{\varphi_\Omega}(x)\lvert \ \text{ a.e  in }\mathbb{R}^N.$$
\end{lemma}
\begin{proof}
	Let $\varphi\in W^{s,p(.,.)}(\Omega)$, since $\overline{p}\in C_+(\overline{\Omega})$, Theorem 3.4.12 in \cite{zbMATH05854634} implies that there exists a sequence $\{\varphi_n\}_{n\in\mathbb{N}}\subset C^{\infty}_{0}(\Omega)$ such that $\varphi_n\longrightarrow\varphi$ strongly in $L^{\overline{p}(.)}(\Omega)$ and $$\lvert \varphi_n(x)\lvert \leqslant\lvert \varphi_{n+1}(x)\lvert \ a.e\ in\ \Omega \text{ for all  }\ n\in\mathbb{N} .$$
	Thus there exists a subsequence still denoted by $\{\varphi_n\}_{n\in\mathbb{N}}$ such that 
	$$\varphi_n\longrightarrow\varphi\ a.e\ in\ \Omega.$$
	Then $$\widetilde{\big(\varphi_n\big)}_{\Omega}\longrightarrow\widetilde{\varphi}_{\Omega}\ a.e\ in\ \mathbb{R}^N.$$
	By Lemma \ref{tracebergot} and the fact that $\varphi_n$ has compact support in $\Omega$, we get 
	$$\widetilde{\big(\varphi_n\big)}_{\Omega}\in W^{s,p(.,.)}_0(\Omega)\ \text{ for all  } n\in\mathbb{N},$$
and the proof is complete.
\end{proof}

\section{The fractional $p(x,.)$- Poisson problem \eqref{poi}}\label{sec3}

In this section, we are concerned with the existence and uniqueness of the weak solution to the fractional $p(x,.)$-Poisson problem.
More precisely, let $p\in \mathfrak{C}(\mathbb{R}^N)$  and $r\in C_+(\overline{ \Omega})$, such that 
\begin{equation}\label{rexp}
	 r(x)<p^*_s(x) \text{ in } \overline{ \Omega},
\end{equation}for $h\in L^{r^{\prime}(.)}(\Omega)$ and $g\in W^{s,p(.,.)}(\mathbb{R}^N) $, we  consider the following fractional nonhomogeneous Poisson equation :
\begin{equation}\tag{$P_{h,g}$}
\left\{\begin{aligned}
\Big(-\Delta_{p(x,.)}\Big)^s u (x)&=h(x)  &\text { in }& \Omega, \\
u &=g &\text { in }& \mathbb{R}^N \setminus\Omega.
\end{aligned}\right.
\end{equation}
We denote by  $\mathcal{K}_{g}(\Omega)=\{u\in W^{s,p(.,.)}(\mathbb{R}^N):\  u-g\in W^{s,p(.,.)}_0(\Omega) \}$, i.e., the class where we are seeking solutions to the Problem \eqref{poi}. To define a weak solution to the Problem \eqref{poi} we denote by $\mathcal{L}$ the operator associated to the fractional $p(x,.)$-Laplacian $(-\Delta_{p(x,.)}^s)$ defined as

\begin{center}
	$\displaystyle \begin{array}{ l c c c c c c }
	\ \ \ \ \mathcal{L} :\  & W^{s,p(.,.)}(\mathbb{R}^N) & \longrightarrow  & \big(W^{s,p(.,.)}(\mathbb{R}^N)\big)^{*} &  &  & \\
	& u & \longrightarrow  & \mathcal{L}( u) : & W^{s,p(.,.)}(\mathbb{R}^N) & \longrightarrow  & \mathbb{R}\\
	&  &  &  & \varphi  & \longrightarrow  & \langle L( u) ,\varphi \rangle, 
	\end{array}$
\end{center}
such that 
$$\langle \mathcal{L}( u) ,\varphi \rangle =\displaystyle\int_{\mathbb{R}^N}\int_{\mathbb{R}^N}\frac{\lvert u(x)-u(y)\lvert ^{p(x,y)-2}(u(x)-u(y))(\varphi(x)-\varphi(y))}{\lvert x-y\lvert^{N+sp(x,y)}} dxdy.$$
Where $\big(W^{s,p(.,.)}(\mathbb{R}^N)\big)^{*}$ is the dual space of $W^{s,p(.,.)}(\mathbb{R}^N).$
\begin{definition}
	We say that $u\in\mathcal{K}_g(\Omega)$ is a weak solution of Problem \eqref{poi}, if 
	$$  \langle \mathcal{L}( u) ,\varphi \rangle=\displaystyle\int_{\Omega}h(x)\varphi(x) dx, $$ for all $\varphi \in W^{s,p(.,.)}_0(\Omega) .$
\end{definition}
The next lemma  will be useful in the subsequent proofs.
\begin{proposition}[\cite{bahrob}]\label{p1}
	Let  $u\in W^{s,p(.,.)}(\mathbb{R}^N)$ and $v\in W^{s,p(.,.)}_0(\Omega)$. Then $$\langle \mathcal{L}( u) ,v \rangle = 2 \int_\Omega v(x) (-\Delta_{p(x,.)}^s)u(x)dx.$$
\end{proposition}
We shall use a variational argument to prove the existence of the weak solutions for Problem
\eqref{poi}. We consider in $W^{s,p(.,.)}(\mathbb{R}^N)$ the energy functional associated with Problem \eqref{poi}.

\begin{equation}\label{F}\tag{$\mathcal{F}$}
\mathcal{F}(u)=\displaystyle\int_{\mathbb{R}^N}\int_{\mathbb{R}^N}\frac{\lvert u(x)-u(y)\lvert ^{p(x,y)}}{p(x,y)\lvert x-y\lvert ^{N+sp(x,y)}} dxdy-\int_{\Omega}h(x)u(x)dx. 
\end{equation}

We recall the definition of minimizer of the functional \ref{F}.
\begin{definition}
	Let $g\in W^{s,p(.,.)}(\mathbb{R}^N) $. A function $u\in\mathcal{K}_{g}(\Omega)$ is a minimizer of the functional $\mathcal{F}$ over $\mathcal{K}_{g}(\Omega)$ if
	\begin{center}
		$\mathcal{F}(u)\leqslant \mathcal{F}(u+\varphi)$ for all $\varphi \in W^{s,p(.,.)}_0(\Omega)$.
		
	\end{center}
\end{definition}
It is standard to show, that a minimizer  of \ref{F} over the class of functions $\mathcal{K}_{g}(\Omega)$ is a weak solution to Problem  \eqref{poi} and vice versa, as stated in Theorem \ref{Ex} below.
\vspace{0.5cm}
\begin{theorem}\label{Ex}
	 Let $\Omega$ be a bounded smooth domain, $ p\in \mathfrak{C}(\mathbb{R}^N)$ with $sp^+<N$ and let $g\in W^{s,p(.,.)}(\mathbb{R}^N) $.  Assume that $r\in C_+(\Omega)$ satisfying \eqref{rexp}. Then for all  $h\in L^{r^{\prime}(.)}(\Omega)$, there exists a unique minimizer $u$ of \ref{F}. Moreover, a function $u\in \mathcal{K}_{g}(\Omega)$ is a minimizer of \ref{F} over $\mathcal{K}_{g}(\Omega)$ if and only if it is a weak solution to Problem \eqref{poi}.
\end{theorem}
\vspace{0.5cm}
In order to prove Theorem \ref{Ex}, we need the following lemma.
\begin{lemma}\label{lemma1}
Let $\Omega$ be a bounded smooth domain, $ p\in \mathfrak{C}(\mathbb{R}^N)$ with $sp^+<N$ and  $ u\in \mathcal{K}_{g}(\Omega) $. Then, there exist positive constants $C_0,\ C_1$ and $C_2$ such that the following inequalities hold:
	
	\begin{enumerate}
		\item $\|u\|_{\overline{p}(.),\Omega}   \leqslant \|u\|_{\overline{p}(.),\mathbb{R}^N}\leqslant C_0\bigg( \|u\|_{\overline{p}(.),\Omega} +\|g\|_{\overline{p}(.),\mathbb{R}^N\setminus \Omega}            \bigg) $ .\label{firstinequality1}
		
		\item 	$\|u\|_{\overline{p}(.),\Omega}\leqslant C_1 [u]_{s,p(.,.),\mathbb{R}^N}$.\label{firstinequality2}
		\item   $\|u\|_{s,p(.,.),\mathbb{R}^N} \leqslant C_2 \big( [u]_{s,p(.,.),\mathbb{R}^N}+\|g\|_{\overline{p}(.),\mathbb{R}^N\setminus\Omega}\big).$ \label{firstinequality3}

	\end{enumerate}

\end{lemma}

\begin{proof}
	Let $ u\in \mathcal{K}_{g}(\Omega) $. The first inequality in \eqref{firstinequality1} is obvious.
	Using the fact that the Orlicz norm $\|\lvert \cdot\|\lvert _{L^{\overline{p}(.)}}$ and the Luxembourg norm $\|\cdot\|_{\overline{p}(.)}$ are equivalent \cite[Theorem 4.8.5]{alois}, we have:
	\begin{eqnarray*}
		\|\lvert u\lvert \|_{L^{\overline{p}(.)}(\mathbb{R}^N)} &=& \displaystyle\sup_{\|\varphi\|_{L^{\bar{p}^{'}(.)}(\mathbb{R}^N)}\leqslant 1}\bigg\lvert \int_{\mathbb{R}^N} u\varphi dx\bigg\lvert ,\\
		&\leqslant& \displaystyle\sup_{\rho_{\overline{p}^{\prime}(.),\mathbb{R}^N}(\varphi)\leqslant 1}\bigg( \int_{\Omega} \lvert u\varphi\lvert  dx +\int_{\mathbb{R}^N\setminus \Omega} \lvert g\varphi\lvert  dx\bigg),\\
		&\leqslant& 2 \displaystyle \sup_{\rho_{\overline{p}^{\prime}(.),\mathbb{R}^N}(\varphi)\leqslant 1} \bigg(  \|u \|_{\overline{p}(.),\Omega}  \|\varphi \|_{q^{\prime}(.),\Omega} + \|g \|_{\overline{p}(.),\mathbb{R}^N\setminus\Omega} \|\varphi \|_{q^{\prime}(.),\mathbb{R}^N\setminus\Omega} \bigg) ,\\
		&\leqslant& C_0\bigg(\|u\|_{\overline{p}(.),\Omega}+\|g\|_{\overline{p}(.),\mathbb{R}^N\setminus\Omega} \bigg).
\end{eqnarray*}
The first assertion is proved.

For \eqref{firstinequality2}, we have to distinguish two cases depending on $g$. In the case $\|g\|_{\overline{p}(.),\mathbb{R}^N}=0$, inequality \eqref{firstinequality2} is an immediate consequence of Lemma \ref{poincaré}. \\
On the other hand, in the case when $\|g\|_{\overline{p}(.),\mathbb{R}^N}\neq 0$,
we set  $$\mathcal{A}=\{ u\in \mathcal{K}_{g}(\Omega),\ \|u\|_{\overline{p}(.),\Omega}=\|g\|_{\overline{p}(.),\Omega}\}.$$
So it suffices to prove that $\displaystyle\inf_{u\in\mathcal{A}}[u]_{s,p(.,.),\mathbb{R}^N} $ is attained for some $u_0\in\mathcal{A}$.
Taking $\{u_n\}$ a minimizing sequence for $\mathcal{A}$, consequently $\{u_n\}$ is bounded in $L^{\overline{p}(.)}(\Omega)$ and in $W^{s,p(.,.)}(\mathbb{R}^N)$, hence $\{u_n\}$ is bounded in $W^{s,p(.,.)}(\Omega)$.

In view of Theorem \ref{emb}, there exists $u_0\in L^{\overline{p}(.)}(\Omega)$ such that, up to subsequence still denoted by $u_n$, $u_n\longrightarrow u_0$ strongly in $L^{\overline{p}(.)}(\Omega)$ and $\|u_0\|_{\overline{p}(.),\Omega}=\|g\|_{\overline{p}(.),\Omega}.$ 
We can extend $u_0$ to $\mathbb{R}^N$ by setting $u_0=g\ in\  \mathbb{R}^N\setminus \Omega,$ this implies that $u_n(x)\longrightarrow u_0(x)$ a.e. in $\mathbb{R}^N$. Using  Fatou's Lemma, we have 
\begin{equation*}
\int_{\mathbb{R}^{N}\times \mathbb{R}^{N}} \dfrac{\lvert u_{0}(x)-u_{0}(y)\lvert ^{p(x,y)}}{\lvert x-y\lvert ^{N+sp(x,y)}}dx dy \leqslant \liminf_{n\to \infty} \int_{\mathbb{R}^{N}\times \mathbb{R}^{N}} \dfrac{\lvert u_{n}(x)-u_{n}(y)\lvert ^{p(x,y)}}{\lvert x-y\lvert ^{N+sp(x,y)}}dx dy ,
\end{equation*}
and $$[u_0]_{s,p(.,.),\mathbb{R}^N}\leqslant \lim_{n\to \infty}[u_n]_{s,p(.,.),\mathbb{R}^N},$$ hence $u_0\in W^{s,p(.,.)}(\mathbb{R}^N)  $ and $u_0\in \mathcal{A}$.
Consequently, we obtain $\displaystyle [u_0]_{s,p(.,.),\mathbb{R}^N}=\inf_{u\in\mathcal{A}}[u]_{s,p(.,.),\mathbb{R}^N} .$ Moreover $[u_0]_{s,p(.,.),\mathbb{R}^N}>0.$\\
Let $\lambda =\displaystyle\inf_{u\in\mathcal{A}}[u]_{s,p(.,.),\mathbb{R}^N}$ and  $u\in\mathcal{K}_g(\Omega)$ such that $\|u\|_{\overline{p}(.),\Omega}\neq 0$, since $\frac{\big(\|g\|_{\overline{p}(.),\Omega}\big)u}{\|u\|_{\overline{p}(.),\Omega}}\in\mathcal{A},$ we get  $$\Big[\frac{\big(\|g\|_{\overline{p}(.),\Omega}\big)u}{\|u\|_{\overline{p}(.),\Omega}}\Big]_{s,p(.,.),\mathbb{R}^N}\geqslant \lambda.$$
Hence $\|u\|_{\overline{p}(.),\Omega}\leqslant \frac{\|g\|_{\overline{p}(.),\Omega}}{\lambda} [u]_{s,p(.,.),\mathbb{R}^N}$, 
The inequality \eqref{firstinequality3} follows directly from the above inequalities.
\end{proof}

The interest of the above lemma is that it allows us to prove the coercivity of the functional $\mathcal{F}$ on $\mathcal{K}_g(\Omega).$

\begin{proof}{\textit{of Theorem} \ref{Ex}.} 
	One can apply the direct method of calculus of variations. Note that the uniqueness is the direct consequence of the strict convexity of the mapping $t\mapsto t^{p(x,y)}$ which is uniform in $(x,y)\in\mathbb{R}^N\times\mathbb{R}^N$, since $p^->1$.
	
	We now prove the existence. Let $u\in \mathcal{K}_g(\Omega)$, by using  Hölder's inequality, Proposition \ref{Pro}, Theorem \ref{emb} and Lemma \ref{lemma1}, we have  
	\begin{eqnarray*}
		\mathcal{F}(u)&=&\displaystyle\int_{\mathbb{R}^N}\int_{\mathbb{R}^N}\frac{\lvert u(x)-u(y)\lvert ^{p(x,y)}}{p(x,y)\lvert x-y\lvert ^{N+sp(x,y)}} dxdy-\int_{\Omega}h(x)u(x)dx\\
		&\geqslant& \frac{1}{p^+}\int_{\mathbb{R}^N}\int_{\mathbb{R}^N}\frac{\lvert u(x)-u(y)\lvert ^{p(x,y)}}{\lvert x-y\lvert ^{N+sp(x,y)}} dxdy - 2\|h\|_{r^{\prime}(.),\Omega}\|u\|_{r(.),\Omega}\\
		&\geqslant& \frac{1}{p^+} \min \big\{[u]^{p+}_{s,p(.,.),\mathbb{R}^N} ;[u]^{p-}_{s,p(.,.),\mathbb{R}^N}\big\}- C\|u\|_{s,p(.,.),\Omega}\\
		&\geqslant& C^{'} \min\big\{\big(\|u\|_{s,p(.,.)}-\|g\|_{\overline{p}(.),\mathbb{R}^N\setminus\Omega)}\big)^{p^-} ; \big(\|u\|_{s,p(.,.)}-\|g\|_{\overline{p}(.),\mathbb{R}^N\setminus\Omega)}\big)^{p^+}\big\}\\
		& & -C\|u\|_{s,p(.,.),\mathbb{R}^N}.
		\end{eqnarray*}
		Since $1<p^-\leqslant p^+$, then $\mathcal{F}$ is coercive on $\mathcal{K}_g(\Omega)$.
	
	Let $\{u_n\}\subset \mathcal{K}_g(\Omega)$ be a minimizing sequence of $\mathcal{F}$, in view of coercivity of $\mathcal{F}$, $\{u_n\}$ is bounded in $W^{s,p(.,.)}(\mathbb{R}^N)$, thus $\{u_n\}$ is bounded in $W^{s,p(.,.)}(\Omega)$, so by Theorem \ref{emb}  there exists a subsequece  still denoted $\{u_n\}$ converging pointwise a.e. in $\Omega$ to a function $u\in L^{\overline{p}(.)}(\Omega)$. We can extend $u$ to $\mathbb{R}^N$ by setting $u=g$ in $\mathbb{R}^N\setminus \Omega$. Thus by Fatou's Lemma, we have $u\in \mathcal{K}_g(\Omega)$ and 
	\begin{eqnarray*}
		\mathcal{F}(u)&=&\displaystyle\int_{\mathbb{R}^N}\int_{\mathbb{R}^N}\frac{\lvert u(x)-u(y)\lvert ^{p(x,y)}}{p(x,y)\lvert x-y\lvert ^{N+sp(x,y)}} dxdy-\int_{\Omega}h(x)u(x)dx\\
		&\leqslant&\liminf_{n\to\infty} \int_{\mathbb{R}^N}\int_{\mathbb{R}^N}\frac{\lvert u_n(x)-u_n(y)\lvert ^{p(x,y)}}{p(x,y)\lvert x-y\lvert ^{N+sp(x,y)}} dxdy-\int_{\Omega}h(x)u_n(x)dx\\
		&\leqslant& \liminf_{n\to\infty}\mathcal{F}(u_n).
	\end{eqnarray*}
	We deduce that $u$ is actually a minimizer of $\mathcal{F}$ over $\mathcal{K}_g(\Omega)$. Furthermore the fact that $u$ solves the corresponding Euler-Lagrange equation follows by perturbing $u\in\mathcal{K}_g(\Omega)$ with a test function in a standard way. Indeed, assume that $u\in\mathcal{K}_g(\Omega)$ a minimizer of \ref{F} over $\mathcal{K}_g(\Omega)$, for any $v\in W^{s,p(.,.)}_0(\Omega)$, we have 
	
	{\footnotesize \begin{eqnarray*}
			\frac{d}{dt}\mathcal{F}(u+tv)\Big\lvert _{t=0}&=&\displaystyle\int_{\mathbb{R}^N}\int_{\mathbb{R}^N}\frac{1}{p(x,y)}\frac{d}{dt}\frac{\lvert u(x)-u(y)-t(v(x)-v(y))\lvert ^{p(x,y)}}{\lvert x-y\lvert ^{N+sp(x,y)}} dxdy-\int_{\Omega}h(x)v(x)dx\Big\lvert _{t=0}\\
			&=& \langle L( u) ,v \rangle-\int_{\Omega}h(x)v(x)dx,
	\end{eqnarray*}}
	Since $u$ is a minimizer, the term on the left is zero and hence $u\in\mathcal{K}_g(\Omega)$ is a weak solution to Problem \eqref{poi}.
	
	Conversely, let $u\in\mathcal{K}_g(\Omega)$ be a weak solution to Problem \eqref{poi} and $v\in\mathcal{K}_g(\Omega)$. Let $\varphi=u-v\in W^{s,p(.,.)}_0(\Omega)$ , we get 
{\footnotesize 	\begin{eqnarray*}
		0 &=&\langle \mathcal{L}( u) ,\varphi \rangle-\int_{\Omega}h(x)\varphi(x)dx\\
		&=& \displaystyle\int_{\mathbb{R}^{2N}}\frac{\lvert u(x)-u(y)\lvert ^{p(x,y)}}{\lvert x-y\lvert ^{N+sp(x,y)}} dxdy-\int_{\mathbb{R}^{2N}}\frac{\lvert u(x)-u(y)\lvert ^{p(x,y)-2}(u(x)-u(y))(v(x)-v(y))}{\lvert x-y\lvert ^{N+sp(x,y)}} dxdy\\
		& &- \int_{\Omega}h(x)u(x)dx+\int_{\Omega}h(x)v(x)dx\\
		&\geqslant& \displaystyle\int_{\mathbb{R}^N\times\mathbb{R}^N}\frac{\lvert u(x)-u(y)\lvert ^{p(x,y)}}{\lvert x-y\lvert ^{N+sp(x,y)}} dxdy-\displaystyle\int_{\mathbb{R}^N\times\mathbb{R}^N}\frac{p(x,y)-1}{p(x,y)}\frac{\lvert u(x)-u(y)\lvert ^{p(x,y)}}{\lvert x-y\lvert ^{N+sp(x,y)}} dxdy\\
		& & -\displaystyle\int_{\mathbb{R}^N\times\mathbb{R}^N}\frac{1}{p(x,y)}\frac{\lvert v(x)-v(y)\lvert ^{p(x,y)}}{\lvert x-y\lvert ^{N+sp(x,y)}} dxdy-\int_{\Omega}h(x)u(x)dx+\int_{\Omega}h(x)v(x)dx\\
		&\geqslant& \displaystyle\int_{\mathbb{R}^N\times\mathbb{R}^N}\frac{1}{p(x,y)}\frac{\lvert u(x)-u(y)\lvert ^{p(x,y)}}{\lvert x-y\lvert ^{N+sp(x,y)}} dxdy-\int_{\Omega}h(x)u(x)dx\\ &  &-\Bigg(\displaystyle\int_{\mathbb{R}^N\times\mathbb{R}^N}\frac{1}{p(x,y)}\frac{\lvert v(x)-v(y)\lvert ^{p(x,y)}}{\lvert x-y\lvert ^{N+sp(x,y)}} dxdy-\int_{\Omega}h(x)v(x)dx \Bigg)\\
		&\geqslant& \mathcal{F}(u)-\mathcal{F}(v),
	\end{eqnarray*}}
	and hence $u$ is a minimizer of \ref{F} over $\mathcal{K}_g(\Omega)$.
\end{proof}

The weak solution of Problem \eqref{poi} enjoys the $L^{r(.)}-$estimates given in the following theorem.
\begin{theorem}\label{Re}
		Let $p$ and $r$ be as above, and 
	let  $h\in L^{r^{'}(.)}(\Omega)$ and  $g\in W^{s,p(.,.)}(\mathbb{R}^N) $.  Then the unique solution of \eqref{poi}  satisfies  $$\|u\|_{r(.),\Omega}^{p(.,.)-1}\leqslant K_1+K_2\bigg(\|h\|_{r^{\prime}(.),\Omega}\bigg)^{\frac{p^+-1}{p^--1}},$$
	where $K_1=K_1(g)$ and $K_2=K_2(g)$ are positive constants independent of $h$.
\end{theorem}
\begin{proof}For the sake of simplicity, we will refer to two positive constants, $K_1$ and $K_2$, whose values are irrelevant and may change from one line to the next.\\
	Let $u\in\mathcal{K}_g(\Omega)$ be a weak solution to Problem \eqref{poi}. Using Theorem \ref{emb} and Lemma \ref{lemma1}, there exists $C>0$ such that 
	\begin{equation}\label{inegalit2}
		\|u\|_{r(.),\Omega}\leqslant C \bigg([u]_{s,p(.,.),\mathbb{R}^N}+\|g\|_{\overline{p}(.),\mathbb{R}^N\setminus\Omega}\bigg). 
	\end{equation}
	Moreover, since $p^+-1\geqslant p^--1>0$ , then by  \eqref{inegalit2}, there exists some positive  constants $K_1, K_2$ such that
	\begin{equation}\label{inegalit3}
	\bigg(\|u\|_{L^{r(.)}(\Omega)}\bigg)^{p^{-}-1} \leqslant K_1\Big([u]_{s,p(.,.),\mathbb{R}^N}\Big)^{p^{-}-1}+ K_2 \Big(\|g\|_{\overline{p}(.),\mathbb{R}^N\setminus\Omega}\Big)^{p^{-}-1}.
	\end{equation}
	and 
		\begin{equation}\label{inegalit4}
	\bigg(\|u\|_{L^{r(.)}(\Omega)}\bigg)^{p^{+}-1} \leqslant K_1\Big([u]_{s,p(.,.),\mathbb{R}^N}\Big)^{p^{+}-1}+ K_2 \Big(\|g\|_{\overline{p}(.),\mathbb{R}^N\setminus\Omega}\Big)^{p^{+}-1}.
	\end{equation}
	Now, we have to distinguish two cases depending $[u]_{s,p(.,.),\mathbb{R}^N}$. 
	 In the case when  $[u]_{s,p(.,.),\mathbb{R}^N}\leqslant 1$, using \eqref{inegalit4}, we get
	  \begin{eqnarray*}
	\bigg(\|u\|_{L^{r^{(.)}}(\Omega)}\bigg)^{p^{+}-1} &\leqslant&  K_1+K_2 \Big(\|g\|_{\overline{p}(.),\mathbb{R}^N\setminus\Omega}\Big)^{p^{+}-1}\\
	 &\leqslant&  K_1+ K_2 \Big(\|g\|_{\overline{p}(.),\mathbb{R}^N\setminus\Omega}\Big)^{p^{+}-1}+\bigg(\|h\|_{r^{\prime}(.),\Omega}\bigg)^{\frac{p^+-1}{p^--1}}\\
	 &\leqslant& K_1+K_2\bigg(\|h\|_{r^{\prime}(.),\Omega}\bigg)^{\frac{p^+-1}{p^--1}},
	\end{eqnarray*}
and similarly, using \eqref{inegalit3}, we get
 \begin{equation}
	\bigg(\|u\|_{L^{r^{(.)}}(\Omega)}\bigg)^{p^{-}-1}
	\leqslant K_1+K_2\bigg(\|h\|_{r^{\prime}(.),\Omega}\bigg)^{\frac{p^+-1}{p^--1}}.
\end{equation}
Then \begin{eqnarray}\label{eq2}
	\bigg(\|u\|_{L^{r^{(.)}}(\Omega)}\bigg)^{p(.,.)-1} &\leqslant &\max \Bigg\{\bigg(\|u\|_{L^{r^{(.)}}(\Omega)}\bigg)^{p^{-}-1};\bigg(\|u\|_{L^{r^{(.)}}(\Omega)}\bigg)^{p^{+}-1} \Bigg\} \nonumber\\
&\leqslant& K_1+K_2\bigg(\|h\|_{r^{\prime}(.),\Omega}\bigg)^{\frac{p^+-1}{p^--1}}.
\end{eqnarray}

	On the other hand, if  $[u]_{s,p(.,.),\mathbb{R}^N}\geqslant 1$, since $u$ is a minimizer of \ref{F} over $\mathcal{K}_g(\Omega)$, we have
	$$\mathcal{F}(u)\leqslant\mathcal{F}(g),$$
	which implies that 
	\begin{eqnarray*}
		\displaystyle\int_{\mathbb{R}^N\times\mathbb{R}^N}\frac{1}{p(x,y)}\frac{\lvert u(x)-u(y)\lvert ^{p(x,y)}}{\lvert x-y\lvert ^{N+sp(x,y)}} dxdy&\leqslant& \displaystyle\int_{\mathbb{R}^N\times\mathbb{R}^N}\frac{1}{p(x,y)}\frac{\lvert g(x)-g(y)\lvert ^{p(x,y)}}{\lvert x-y\lvert ^{N+sp(x,y)}} dxdy\\ & &-\int_{\Omega}h(x)g(x)dx+\int_{\Omega}h(x)u(x)dx.\\
	\end{eqnarray*}
	Hence $$\displaystyle\rho_{s,p(.,.),\mathbb{R^N}}(u)\leqslant \dfrac{p^{+}}{p^{-}}\rho_{s,p(.,.),\mathbb{R^N}}(g) -p^+\int_{\Omega}h(x)g(x)dx+p^+\int_{\Omega}h(x)u(x)dx. $$
	Using Hölder's inequality and \eqref{inegalit2}, we obtain
	{\footnotesize\begin{eqnarray}
		\displaystyle\rho_{s,p(.,.),\mathbb{R^N}}(u)&\leqslant& \dfrac{p^{+}}{p^{-}}\rho_{s,p(.,.),\mathbb{R^N}}(g)+2p^+\|h\|_{r^{\prime}(.),\Omega} \|u\|_{r(.),\Omega}+2p^+\|h\|_{r^{\prime}(.),\Omega} \|g\|_{r(.),\Omega}\nonumber\\
		&\leqslant& \dfrac{p^{+}}{p^{-}}\rho_{s,p(.,.),\mathbb{R^N}}(g)+C\|h\|_{r^{\prime}(.),\Omega}[u]_{s,p,\mathbb{R}^N} +\|h\|_{r^{\prime}(.),\Omega}\big(C\|g\|_{\overline{p}(.),\mathbb{R}^N\setminus\Omega}+\|g\|_{r(.),\Omega}\big) .  \label{inegalit}
	\end{eqnarray}
}

Since $[u]_{s,p(.,.),\mathbb{R}^N}\geqslant 1$, then, invoking Proposition \ref{Pro} and \eqref{inegalit}, we have
	{\footnotesize $$\Big([u]_{s,p(.,.),\mathbb{R}^N}\Big)^{p^{-}}\leqslant\dfrac{p^{+}}{p^{-}}\rho_{s,p(.,.),\mathbb{R^N}}(g)+C\|h\|_{r^{\prime}(.),\Omega}[u]_{s,p(.,.),\mathbb{R}^N} +\|h\|_{r^{\prime}(.),\Omega}\bigg(C\|g\|_{\overline{p}(.),\mathbb{R}^N\setminus\Omega}+\|g\|_{r(.),\Omega}\bigg).$$}
	Hence {\footnotesize\begin{eqnarray}
		\Big([u]_{s,p(.,.),\mathbb{R}^N}\Big)^{p^{-}-1}&\leqslant& \dfrac{p^{+}}{p^{-}}\rho_{s,p(.,.),\mathbb{R^N}}(g)+\bigg(C + C\|g\|_{\overline{p}(.),\mathbb{R}^N\setminus\Omega}+\|g\|_{L^{r(.)}(\Omega)}\bigg)\|h\|_{r^{\prime}(.),\Omega}\\
		&\leqslant& K_1+K_2 \|h\|_{r^{\prime}(.),\Omega}.
	\end{eqnarray}}
	Since $\frac{p^+-1}{p^--1}>0$, we get
	\begin{equation}\label{ineg4}
	\Big([u]_{s,p(.,.),\mathbb{R}^N}\Big)^{p^{+}-1}\leqslant K_1+K_2\bigg(\|h\|_{r^{\prime}(.),\Omega}\bigg)^{\frac{p^+-1}{p^--1}}.
	\end{equation}
	Combining \eqref{inegalit3},\eqref{inegalit4} and \eqref{ineg4}, we obtain 
	\begin{equation}\label{eq11}
	\bigg(\|u\|_{L^{r^{(.)}}(\Omega)}\bigg)^{p(.,.)-1} \leqslant  K_1+K_2\bigg(\|h\|_{r^{\prime}(.),\Omega}\bigg)^{\frac{p^+-1}{p^--1}}.
	\end{equation}

	Finally combining \eqref{eq11} and \eqref{eq2}, we obtain the desired result.
\end{proof}

Let $h\in L^{r^{\prime}(.)}(\Omega) $ , we denote by $\mathcal{T}(h):=u_{h,g}^\Omega$ the unique solution of Problem  \eqref{poi}. As a consequence of Theorem \ref{emb} and using the same approach as in \cite[Lemma 4.4]{zbMATH06813457}, we get the following property of the operator $\mathcal{T}$:
\begin{lemma}\label{CC}
	Let $p\in \mathfrak{C}(\mathbb{R}^N)$. Then for any $r\in C_+(\Omega)$ satisfying \eqref{rexp} , we have 
	\begin{enumerate}
		\item The mapping $\mathcal{T}:L^{r^{'}(.)}(\Omega)\rightarrow W^{s,p(.,.)}(\mathbb{R}^N)$ is continuous.
		\item The mapping $\mathcal{T}:L^{r^{'}(.)}(\Omega)\rightarrow L^{r(.)}(\Omega)$ is completely  continuous, i.e. $h_n\rightharpoonup h$ in $L^{r^{'}(.)}(\Omega)$ implies $\mathcal{T}(h_n)\longrightarrow\mathcal{T}(h)$ in $L^{r(.)}(\Omega)$.
	\end{enumerate} 
\end{lemma}
\section{Main results and proofs}\label{sec4}

In this section, we are concerned with the study  of the Problem \eqref{P}:

\begin{equation*}
\left\{\begin{aligned}
(-\Delta_{p(x,.)})^s u &=f(x, u)  &\text { on }& \Omega, \\
u &=g &\text { in }& \mathbb{R}^N \setminus\Omega.
\end{aligned}\right.
\end{equation*}
	Where $\Omega\subset \mathbb{R}^N$ be a smooth bounded domain, $s\in(0,1)$. From now on, we assume that 
 $p\in \mathfrak{C}(\mathbb{R}^N)$, with $sp^+<N$.

We make the following assumption:
\begin{itemize}
	\item[($f_1$):] $f: \Omega\times\mathbb{R}\longrightarrow \mathbb{R}$ is a continuous function such that $\lvert f(x;\xi)\lvert \leqslant a(x)+C\lvert \xi\lvert ^{\bar{p}(x)-1}$ for all  $x\in\Omega$ and for all $\xi\in\mathbb{R}^N$, where $C>0$ and $a: \Omega\longrightarrow \mathbb{R}^+,$ such that $a\in L^{\bar{p}^{'}(.)}(\Omega)$. 
\end{itemize}

We define the Nemytsky operator $\mathcal{N}_f$ generated by  $f$, acting on
the measurable function $u:\Omega\longrightarrow \mathbb{R}$ by $$\mathcal{N}_f\big(u\big)(x)=f\big(x,u(x)\big),$$
for all $x\in\Omega$. We say that $u\in W^{s,p(.,.)}(\mathbb{R}^N)$ is a weak solution of Problem \eqref{P}, if 
\begin{equation*}
\left\{\begin{aligned}
u&\in \mathcal{K}_g&\\
\mathcal{N}_f\big(u\big) &\in L^{r^{'}(.)}(\Omega) &\\
\langle L( u) ,\varphi \rangle &=\displaystyle\langle\mathcal{N}_f\big(u\big),\varphi\rangle,&  
\end{aligned}\right.
\end{equation*}
for all $\varphi\in W^{s,p(.,.)}_0(\Omega)$, where $\langle\mathcal{N}_f\big(u\big),\varphi\rangle=\displaystyle\int_\Omega\mathcal{N}_f\big(u\big)(x)\varphi(x)dx$ and
\begin{equation*}\label{Cond}\tag{$R_1$}
\begin{split}
\left\{\begin{aligned}
	r&\in C_+(\overline{ \Omega})&\\
\overline{p}(x)\leqslant \overline{p}^+<&r^-\leqslant r(x)<p^*_s(x)& \forall x\in\overline{ \Omega}.\\
\end{aligned}\right.
\end{split}
\end{equation*}
\begin{lemma}\label{est}
Assume that the assumptions  $(f_1)$ and \eqref{Cond}  hold. Then, for all $u\in W^{s,p(.,.)}(\mathbb{R}^N)$, we have $\mathcal{N}_f\big(u\big) \in L^{r^{'}(.)}(\Omega)$.  Moreover there exist a positive constant $C$ such that $$\|\mathcal{N}_f\big(u\big)\|_{L^{r^{'}(.)}(\Omega)}\leqslant C\bigg(\|a\|_{L^{\bar{p}(.)}(\Omega)}+\big(C^{\gamma}_\Omega\|u\|_{L^{r(.)}(\Omega)}\big)^{p^{+}-1}+\big(C^{\gamma}_\Omega\|u\|_{L^{r(.)}(\Omega)}\big)^{p^{-}-1}\bigg),$$ 
where $C^\gamma_\Omega:= 2\max\big\{\lvert \Omega\lvert ^{\frac{1}{\gamma^+}}; \lvert \Omega\lvert ^{\frac{1}{\gamma^-}}\big\}$ and $\gamma(x):=\frac{\bar{p}(x).r(x)}{r(x)-\bar{p}(x)}$.
\end{lemma}

\begin{proof}
	Let $u\in W^{s,p(.,.)}(\Omega)$, we have 
	
	\begin{eqnarray*}
		\lvert \langle \mathcal{N}_f\big(u\big), v \rangle\lvert &\leqslant& \langle \lvert \mathcal{N}_f\big(u\big)\lvert , \lvert v\lvert  \rangle\\
		&\leqslant& \langle \lvert a\lvert   ,\lvert v\lvert  \rangle +\langle \lvert u\lvert ^{\bar{p}(.)-1} ,\lvert v\lvert  \rangle,
	\end{eqnarray*}
	for all  $v\in L^{r(.)}(\Omega).$ Using Hölder's inequality and Lemma \ref{inclusion}, we have
	
	\begin{eqnarray*}
		\lvert \langle \mathcal{N}_f\big(u\big), v \rangle\lvert &\leqslant& 2\bigg(\|a\|_{L^{\bar{p}^{'}(.)}(\Omega)}\|v\|_{L^{\bar{p}(.)}(\Omega)}+\|\lvert u\lvert ^{\bar{p}(.)-1}\|_{L^{\overline{p}^{'}(.)}(\Omega)}\|v\|_{L^{\bar{p}(.)}(\Omega)}\bigg)\\
		&\leqslant& C_1\bigg(\|a\|_{L^{\bar{p}^{'}(.)}(\Omega)}\|v\|_{L^{r(.)}(\Omega)}+\|\lvert u\lvert ^{\bar{p}(.)-1}\|_{L^{\bar{p}^{'}(.)}(\Omega)}\|v\|_{L^{r(.)}(\Omega)}\bigg).\\
	\end{eqnarray*}
	Thus {\footnotesize
		\begin{eqnarray*}
			\sup_{\|v\|_{L^{r(.)}}\leqslant 1}	\lvert \langle \mathcal{N}_f\big(u\big), v \rangle\lvert 
			&\leqslant& C_1\bigg(\|a\|_{L^{\bar{p}^{'}(.)}}\sup_{\|v\|_{L^{r(.)}}\leqslant 1}\|v\|_{L^{r(.)}}+\|\lvert u\lvert ^{\bar{p}(.)-1}\|_{L^{\bar{p}^{'}(.)}}\sup_{\|v\|_{L^{r(.)}}\leqslant 1}\|v\|_{L^{r(.)}}\bigg)\\
			&\leqslant& C_1\bigg(\|a\|_{L^{\bar{p}^{'}(.)}(\Omega)}+\|\lvert u\lvert ^{\bar{p}(.)-1}\|_{L^{\bar{p}^{'}(.)}(\Omega)}\bigg).
	\end{eqnarray*}}
	By Lemma \ref{edm}, we have  
	\begin{eqnarray*}
		\|\lvert u\lvert ^{\bar{p}(.)-1}\|_{L^{\bar{p}^{'}(.)}(\Omega)}&\leqslant& \max\bigg\{\|u\|^{\bar{p}^+-1}_{L^{\bar{p}^{'}(\bar{p}-1)}(\Omega)};\|u\|^{\bar{p}^--1}_{L^{\bar{p}^{'}(\bar{p}-1)}(\Omega)}\bigg\}\\
		&\leqslant& \|u\|^{\bar{p}^+-1}_{L^{\bar{p}^{'}(\bar{p}-1)}(\Omega)}+\|u\|^{\bar{p}^--1}_{L^{\bar{p}^{'}(\bar{p}-1)}(\Omega)}\\
		&\leqslant& \|u\|^{\bar{p}^+-1}_{L^{\bar{p}(.)}(\Omega)}+\|u\|^{\bar{p}^--1}_{L^{\bar{p}(.)}(\Omega)}.
	\end{eqnarray*}
	Hence 
	\begin{eqnarray*}
		\sup_{\|v\|_{L^{r(.)}(\Omega)}\leqslant 1}	\lvert \langle \mathcal{N}_f\big(u\big), v \rangle\lvert 
		&\leqslant& C_1\bigg(\|a\|_{L^{\bar{p}^{'}(.)}(\Omega)}+\|u\|^{\bar{p}^+-1}_{L^{\bar{p}(.)}(\Omega)}+\|u\|^{\bar{p}^--1}_{L^{\bar{p}(.)}(\Omega)}\bigg).
	\end{eqnarray*}
   Since $$\frac{1}{\gamma(x)}+\frac{1}{r(x)}=\frac{1}{\overline{p}(x)}\ \text{ for all  }\ x\in \overline{ \Omega},$$
	then, using Lemmas \ref{holder2} and \ref{Cara},  and the fact that the Luxembourg norm and the Orlicz norm are equivalent, we have 
{\footnotesize	\begin{eqnarray*}
		\|\mathcal{N}_f\big(u\big)\|_{L^{r^{'}(.)}(\Omega)}&\leqslant& C\bigg(\|a\|_{L^{\bar{p}^{'}(.)}(\Omega)}+\|u\|^{\bar{p}^+-1}_{L^{\bar{p}(.)}(\Omega)}+\|u\|^{\bar{p}^--1}_{L^{\bar{p}(.)}(\Omega)}\bigg)\\
		&\leqslant& C\bigg(\|a\|_{L^{\bar{p}^{'}(.)}(\Omega)}+\big(\|1\|_{L^\gamma(\Omega)}\|u\|_{L^{r(.)}(\Omega)}\big)^{p^{+}-1}+\big(\|1\|_{L^\gamma(\Omega)}\|u\|_{L^{r(.)}(\Omega)}\big)^{p^{-}-1}\bigg)\\
		&\leqslant&C\bigg(\|a\|_{L^{\bar{p}^{'}(.)}(\Omega)}+\big(C^\gamma_\Omega\|u\|_{L^{r(.)}(\Omega)}\big)^{p^{+}-1}+\big(C^\gamma_\Omega\|u\|_{L^{r(.)}(\Omega)}\big)^{p^{-}-1}\bigg),
	\end{eqnarray*}}
	and the proof is complete.
\end{proof}

\begin{lemma}\label{CB}
	Assume that the assumptions $(f_1)$ and (\ref{Cond}) hold. Then the operator 
	$$\mathcal{N}_f:L^{\overline{p}(.)}(\Omega)\longrightarrow L^{r^{'}(.)}(\Omega)$$ is continuous and bounded.
\end{lemma}
\begin{proof} 
	We first notice that
	$$\frac{\overline{p}^{'}(x)-r^{\prime}(x)}{\overline{p}^{'}(x)}+\frac{r^{'}(x)}{\overline{p}^{'}(x)}=1\ \text{ for all  }\ x\in \overline{ \Omega}.$$
	Thus by Young's inequality, we get
	\begin{eqnarray}
		\lvert f(x,t)\lvert &\leqslant& a(x) +C\lvert t\lvert ^{\bar{p}(x)-1}\nonumber\\
		&\leqslant& a(x) +\frac{\overline{p}^{'}(x)-r^{\prime}(x)}{\overline{p}^{'}(x)} C^{\frac{\overline{p}^{'}(x)}{\overline{p}^{'}(x)-r^{\prime}(x)}}+\frac{r^{'}(x)}{\overline{p}^{'}(x)} \lvert t\lvert ^\frac{\overline{p}^{'}(x)}{r^{'}(x)},\nonumber\\
		 &\leqslant & a(x)+C^{'}+\bigg(\frac{r^{'}}{\overline{p}^{'}}\bigg)^+\lvert t\lvert ^\frac{\overline{p}^{'}(x)}{r^{'}(x)},\label{3}
	\end{eqnarray}
	where $C^{'}=\Big(\frac{\overline{p}^{'}-r^{\prime}}{\overline{p}^{'}}\Big)^+\max \Bigg\{C^{\bigg(\frac{\overline{p}^{'}}{\overline{p}^{'}-r^{\prime}}\bigg)^{+}};C^{\bigg(\frac{\overline{p}^{'}}{\overline{p}^{'}-r^{\prime}}\bigg)^{-}}\Bigg\}.$\\
	Since $a\in L^{\bar{p}^{'}}(\Omega)$, then by \eqref{Cond} and using Lemmas \ref{inclusion} and \ref{Cara}, we obtain 
	\begin{eqnarray*}
		\big\|  a	+C^{'}\big\|_{L^{r^{'}}(\Omega)}&\leqslant & \|a\|_{L^{r^{'}}(\Omega)}+C^{'}\|1\|_{L^{r^{'}}(\Omega)}\\
		&\leqslant&\|a\|_{L^{r^{'}}(\Omega)}+C^{'}\max\bigg\{\lvert \Omega\lvert ^{\frac{1}{r^{'+}}}; \lvert \Omega\lvert ^{\frac{1}{r^{'-}}}\bigg\}<\infty
	\end{eqnarray*}
	Then \begin{equation}\label{4}
	\bigg(a(x)+C^{'}\bigg)\in L^{r^{'}}(\Omega).
	\end{equation}
	Combining \eqref{3}, \eqref{4} and Lemma \ref{nemytsky}, we deduce that $\mathcal{N}_f$ is continuous and bounded from $L^{\overline{p}(.)}(\Omega)$ to $L^{r^{'}(.)}(\Omega)$.
\end{proof}

Now, we can state the existence result for  domains with small Lebesgue measure, using Schauder's fixed point theorem.
\begin{theorem}
	Assume the assumptions $(f_1)$ and (\ref{Cond}) hold. If $\lvert \Omega \lvert $ is small enough, 
	then  for every $g\in W^{s,p(.,.)}(\mathbb{R}^N)$, the Problem \eqref{P} has at least one solution  $u^{\Omega}_{g}$.
\end{theorem}

\begin{proof}
	We define $$\mathcal{J}(h):=\mathcal{N}_f\circ \mathcal{T}(h).$$
	Then , by Lemmas \ref{CC} and \ref{CB}, we have $\mathcal{J}:L^{r^{'}}(\Omega)\longrightarrow L^{r^{'}}(\Omega)$ is completely continuous.\\
	In order to prove that Problem \eqref{P} has a weak solution $u_g^\Omega\in\mathcal{K}_g(\Omega)$, it is sufficient to prove that the equation 
	\begin{equation}\label{Fix}
	\mathcal{J}(h)=h
	\end{equation}
	has a solution in $L^{r^{'}}(\Omega)$.	Indeed, if $h\in L^{r^{'}}(\Omega)$ is a solution of (\ref{Fix}), then by Theorem \ref{Ex}, there exist a unique $u_{h,g}^{\Omega}=\mathcal{T}(h)\in\mathcal{K}_g(\Omega)$ such that $$  \langle L( u_{h,g}) ,\varphi \rangle=\displaystyle\int_{\Omega}h\varphi dx, \text{ for all }\varphi \in W^{s,p(.,.)}_0(\Omega).$$
	On the other hand, by the definition of the operator $\mathcal{J}$, we have 
	\begin{eqnarray*}
		\langle L( u_{h,g}^{\Omega}) ,\varphi \rangle&=&\displaystyle\int_{\Omega}h\varphi dx\\
		&=&\displaystyle\int_{\Omega}\mathcal{J}(h)\varphi dx\\
		&=&   \displaystyle\int_{\Omega}\big(\mathcal{N}_f\circ \mathcal{T}(h)\big)\varphi dx\\
		&=&\langle\mathcal{N}_f\big(u_{h,g}^{\Omega}\big),\varphi\rangle,
	\end{eqnarray*}
	for all $\varphi\in W^{s,p(.,.)}_0(\Omega)$. Thus $u_{h,g}^{\Omega}$ is a solution to Problem \eqref{P}. Therefore, proving that the  Problem \eqref{P} has a  solution remains to show that the compact operator $\mathcal{J}$ has a fixed point.\\
	Let us set $$\mathcal{E}_M:=\{h\in L^{r^{'}}(\Omega);\ \|h\|_{L^{r^{'}}(\Omega)}\leqslant M\},$$
	where $M$ is a constant to be chosen later . Notice that $\mathcal{E}_M$ is a non empty closed convex subset of $L^{r^{'}}(\Omega)$. By Lemma \ref{est}, we have
	\begin{eqnarray*}
		\|\mathcal{J}(h)\|_{L^{r^{'}(.)}(\Omega)} &=&	\|\mathcal{N}_f\big(\mathcal{T}(h)\big)\|_{L^{r^{'}(.)}(\Omega)}\\
		&=&	\|\mathcal{N}_f\big(u_{h,g}^{\Omega}\big)\|_{L^{r^{'}(.)}(\Omega)}\\
		&\leqslant& C\bigg(\|a\|_{L^{\bar{p}(.)}(\Omega)}+\big(C_\Omega\|u_{h,g}^{\Omega}\|_{L^{r(.)}(\Omega)}\big)^{p^{+}-1}+\big(C_\Omega\|u_{h,g}^{\Omega}\|_{L^{r(.)}(\Omega)}\big)^{p^{-}-1}\bigg).
	\end{eqnarray*}
	By Theorem \ref{Re}, we obtain 
	$$\|u_{h,g}^{\Omega}\|_{L^{r(.)}(\Omega)}^{p^+-1}\leqslant K_1+ K_2\bigg(\|h\|_{r^{\prime}(.),\Omega}\bigg)^{\frac{p^+-1}{p^--1}} \text{ and }  \|u_{h,g}^{\Omega}\|_{L^{r(.)}(\Omega)}^{p^--1}\leqslant K_1+K_2\bigg(\|h\|_{r^{\prime}(.),\Omega}\bigg)^{\frac{p^+-1}{p^--1}} ,$$
	where $K_1>0$ and $K_2>0$ are independent of $h.$ Hence
	 {\small\begin{eqnarray*}
			\|\mathcal{J}(h)\|_{L^{r^{'}(.)}(\Omega)} &\leqslant&C\bigg(\|a\|_{L^{\bar{p}(.)}(\Omega)}+\big(C_\Omega\big)^{p^{+}-1}\Big(K_1+K_2\big(\|h\|_{r^{\prime}(.),\Omega}\big)^{\frac{p^+-1}{p^--1}}\Big)\\
			& & +\big(C_\Omega\big)^{p^{-}-1}\Big(K_1+K_2\big(\|h\|_{r^{\prime}(.),\Omega}\big)^{\frac{p^+-1}{p^--1}}\Big)\bigg)\\
			&\leqslant& C\Big(\|a\|_{L^{\bar{p}(.)}(\Omega)} +K_1\big(C_\Omega\big)^{p^{+}-1}+K_1\big(C_\Omega\big)^{p^{-}-1}\Big)\\
			& &+M^{\frac{p^+-1}{p^--1}}K_2\bigg(\big(C_\Omega\big)^{p^{+}-1}+\big(C_\Omega\big)^{p^{-}-1}\bigg).
	\end{eqnarray*}}
Where, 
	Since $\lvert \Omega\lvert $ is small enough, we can choose  $$M_\Omega:=\Bigg(\frac{ C\big(\|a\|_{L^{\bar{p}(.)}(\Omega)} +K_1\big(C_\Omega\big)^{p^{+}-1}+K_1\big(C_\Omega\big)^{p^{-}-1}\big)}{1-K_2\big(\big(C_\Omega\big)^{p^{+}-1}+\big(C_\Omega\big)^{p^{-}-1}\big)}\Bigg)^{\frac{p^--1}{p^+-1}}>0.$$
	
	Consequently, $$\mathcal{J}(\mathcal{E}_M)\subset \mathcal{E}_M, \text{ for all  } M\geqslant M_\Omega .$$
	The hypothesis of the Schauder's fixed point theorem are satisfied. Hence there exists at least one solution to the Problem \eqref{P}.
\end{proof}

Now, we are in a position  to state our existence result with no restriction on the measure of $\Omega$. 

\begin{theorem}
  Assume that the assumptions $\big(f_1\big)$ and \eqref{Cond} hold. Then for every $g\in W^{s,p(.,.)}(\mathbb{R}^N)$, there exists at least one solution $u_g^\Omega$ to Problem \eqref{P}.
\end{theorem}

\begin{proof}
	Let $g\in W^{s,p(.,.)}(\mathbb{R}^N)$. We define a sequence of open subsets $\{\Omega_j\}_{j=0}^{j=M+1}$ of $\Omega$ such that $\Omega_0=\Omega,\ \Omega_{M+1}=\varnothing,\ \Omega_{j+1}\subset\Omega_{j} $ and $V_j=\Omega_{j}\setminus\overline{\Omega_{j+1}}$ with  $|V_j|$ being small enough.
	
	 We define the functions $\{w_j\}_{j=0}^{j=M}$ by
	
	\begin{equation*}
	\left\{\begin{aligned}
	w_0&=u_{g}^{V_0}&\\
	w_{j+1}&=u^{V_{j+1}}_{w_j}&\hspace{1cm} (j=0,\ldots,M-1).
	\end{aligned}\right.
	\end{equation*}
	Then $w_{M}\in\mathcal{K}_g(\Omega)$ is a solution of Problem \eqref{P}. Indeed, let $\varphi\in W^{s,p(.,.)}_0(\Omega)$. By Proposition \ref{p1}, we have 
	\begin{eqnarray*}
		\langle L( w_{M}) ,\varphi \rangle&=&  2 \int_\Omega \varphi(x) (-\Delta_{p(x,.)})^s w_{M}(x)dx\\
		&=&   \sum_{j=0}^{j=M}2\int_{V_j} \varphi(x) (-\Delta_{p(x,.)})^s w_{M}(x)dx\\
		&=&\sum_{j=0}^{j=M}2\int_{V_j} \widetilde{\varphi}_{V_j}(x) (-\Delta_{p(x,.)})^s w_{j}(x)dx.
	\end{eqnarray*}
	By Lemma \ref{dx} and the fact that $\varphi\in W^{s,p(.,.)}_0(\Omega)$,  there exists a sequence \begin{equation}\label{eq1}
	\{\varphi_{j,n}\}_{n\in\mathbb{N}}\subset W^{s,p(.,.)}_0(V_j)
	\end{equation} such that $$\lvert \varphi_{j,n}\lvert \nearrow \lvert \widetilde{\varphi_{V_j}}\lvert \ \text{ a.e  in }\mathbb{R}^N,\ for\  j=0,\ldots,M.$$
	Therefore, by Proposition \ref{p1}, the definition of $w_{j}$ and $(\ref{eq1})$, together with the monotone convergence theorem, we have 
	\begin{eqnarray*}
		2\int_{V_j} \widetilde{\varphi}_{V_j}(x) (-\Delta_{p(x,.)})^s w_{j}(x)dx&=& \lim_{n\to+\infty} 2\int_{V_j} \varphi_{j,n}(x) (-\Delta_{p(x,.)})^s w_{j}(x)dx\\
		&=& \lim_{n\to+\infty}	\langle L( w_{j}) , \varphi_{j,n}\rangle\\
		&=& \lim_{n\to+\infty} \int_{V_j}\mathcal{N}_f\big(w_{j}\big)\varphi_{j,n}\\
		&=& \int_{V_j}\mathcal{N}_f\big(w_{j}\big)\widetilde{\varphi}_{V_j}\\
		&=&\int_{V_j}\mathcal{N}_f\big(w_{j}\big)\varphi.
	\end{eqnarray*} 
	Thus 
	\begin{eqnarray*}
			\langle L( w_{M}) ,\varphi \rangle 
		&=&   \sum_{j=0}^{j=M}\int_{V_j}\mathcal{N}_f\big(w_{j}\big)\varphi\\
		&=&   \sum_{j=0}^{j=M}\int_{V_j}\mathcal{N}_f\big(w_{M}\big)\varphi\\
		&=& \int_{\Omega}\mathcal{N}_f\big(w_{M}\big)\varphi.
	\end{eqnarray*}
	Hence $w_{M}$ is a solution to Problem \eqref{P}, and this complete the proof of our theorem.
\end{proof}

\end{document}